# A Family of Quantile Dependence Coefficients

## A. Dastbaravarde [a] and A. Dolati [b]

[a,b] *Department of Statistics, Yazd University, Yazd, Iran.*
adolati@yazd.ac.ir

### Abstract

A popular measure of association is the tail dependence coefficient which measures the strength of dependence in either the lower-left or upper-right tail of a bivariate distribution. In this paper, we develop the idea of quantile dependence, which generalizes the notion of tail dependence and could be used to detect dependence in specific regions of the domain of a joint distribution function. Properties of the proposed quantile dependence coefficient are studied and several examples illustrate our results.

## 1    Introduction

Most of the dependence measures associate the entire distribution of random variables. However, the dependence between the upper tail of the distribution may be different than the lower tail of the distribution. For example, two random variables with weak dependence between mid-range values, but strong dependence in high or low values. The concept of tail dependence refers to the amount of dependence in the upper-right-quadrant tail or lower-left-quadrant tail of a bivariate distribution. The notion of tail dependence coefficient has been introduced by Sibuya [17]. Tail dependence coefficients (upper and lower) are defined as the asymptotic probability that two extremal events occur simultaneously. Lower and upper tail dependence coefficients are important for modeling the dependence structure between random variables in many fields. For example, in finance [7, 9, 19], tail dependence coefficients are used to model the joint tail risk of asset returns. This is important for risk management because extreme events in financial markets can have significant impacts on portfolios. By understanding the dependence structure between assets in the lower and upper tails, investors can better manage their risk exposure. In insurance [6, 10], modeling the joint tail risk of insurance claims is important for insurers because they need to be able to estimate the probability of extreme events that could lead to large payouts. By understanding the dependence structure between different types of insurance claims in the lower and upper tails, insurers can better manage their risk exposure and set appropriate premiums. In environmental science, tail dependence coefficients are used to model the joint tail risk of extreme weather events which is important for understanding the potential impacts of climate change on ecosystems, infrastructure, and human health. By understanding the dependence structure between different types of extreme weather events in the lower and upper tails, scientists can better predict the frequency and severity of extreme events in the future [8, 12]. The most common definition of tail dependence coefficients as provided in [4] is the following.

**Definition 1** *Let $X$ and $Y$ be two continuous random variables with the cumulative distribution functions (cdf) $F$ and $G$, respectively. The upper tail dependence coefficient $\lambda_U$ of $(X, Y)$ is defined by*

$$\lambda_U = \lim_{t \to 1^-} P\{Y > G^{-1}(t) | X > F^{-1}(t)\}, \tag{1.1}$$

*and the lower tail dependence coefficient $\lambda_L$ of $(X, Y)$ is defined by*

$$\lambda_L = \lim_{t \to 0^+} P\{Y \leq G^{-1}(t) | X \leq F^{-1}(t)\}, \tag{1.2}$$

*provided that the above limits exist, where $F^{-1}(u) = \inf\{x \in \ I\!R : F(x) \geq u\}$, is the quantile function or generalized inverse of the cdf $F$.*

The upper tail dependence coefficient indicates the limit of the probability that the random variable $Y$ exceeds a high quantile of its distribution, given that the random variable $X$ exceeds a high quantile of its distribution. A similar interpretation holds for the lower tail dependence coefficient. The upper and lower tail dependence coefficients can be expressed in terms of the copula $C$ of $(X, Y)$ given the well-known Sklar's Theorem [18] via $H(x, y) = C(F(x), G(y))$, by

$$\lambda_L = \lambda_L(C) = \lim_{t \to 0^+} \frac{C(t,t)}{t}, \tag{1.3}$$

and

$$\lambda_U = \lambda_U(C) = \lim_{t \to 0^+} \frac{\overline{C}(1-t, 1-t)}{t} = \lim_{t \to 0^+} \frac{\hat{C}(t,t)}{t} = \lambda_L(\hat{C}), \tag{1.4}$$

where $\overline{C}(u, v) = 1 - u - v + C(u, v)$ is the survival function of $C$, and $\hat{C}(u, v) = \overline{C}(1 - u, 1 - v)$ is the survival copula associated with the copula $C$. As the tail dependence coefficients can be expressed via a copula, many properties of copulas, for example, invariance under strictly increasing transformations of the margins, also apply to the tail dependence coefficient. Copula $C$ has an upper tail dependence if $\lambda_U \in (0, 1]$ and no upper tail dependence if $\lambda_U = 0$. Copula $C$ has lower tail dependence if $\lambda_L \in (0, 1]$ and no lower tail dependence if $\lambda_L = 0$. For the independence copula $\Pi(u, v) = uv$ we have tail independence ($\lambda_L = \lambda_U = 0$) and for the Fréchet-Hoeffding upper bound copula $M(u, v) = \min(u, v)$ [11] we have perfect tail dependence ($\lambda_U = \lambda_L = 1$). Tail dependence coefficients for the most popular families of copulas are available in [3]. For more information about tail dependence coefficients, their properties, applications, and generalizations see, e.g., [1, 5, 12, 13, 15, 16]. Lower and upper tail dependence coefficients evaluate dependence in the tails of between variables. However, one may be interested in analyzing dependence in some specific parts of the distribution, rather than tails. As an example, in financial markets, the dependence between asset prices may be significantly higher during periods of crisis. This dependency breakdown can occur in any part of the distribution domain. The following example illustrates the problem.

**Example 1** *Consider a pair $(X, Y)$ of continuous random variables with $X \sim N(0, 1)$ and*

$$Y|(X = x) \sim N(\beta_0 + \beta_1 x, x^2), \qquad \beta_0, \beta_1 \in R.$$

*Figure 1 shows the scatter plots 10000 sample points generated from a pair $(X, Y)$ distributed as this model, with different values of the parameters. The left panels show the scatterplot of the data points and the right panels show their corresponding normalized ranks. The strong dependence of the variables in the areas other than the tails are specified in the plots. The*

Pearson's correlation coefficient of $(X,Y)$ for this model is given by $\rho_{X,Y} = \frac{\beta_1}{\sqrt{1+\beta_1^2}}$. We note that for the case $\beta_1 = 0$, $\rho_{X,Y} = 0$ but there is strong dependence in the point $(0,0)$.

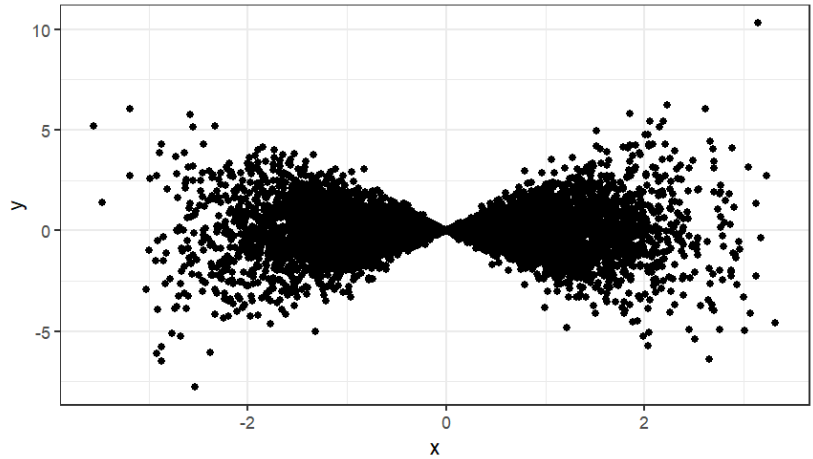

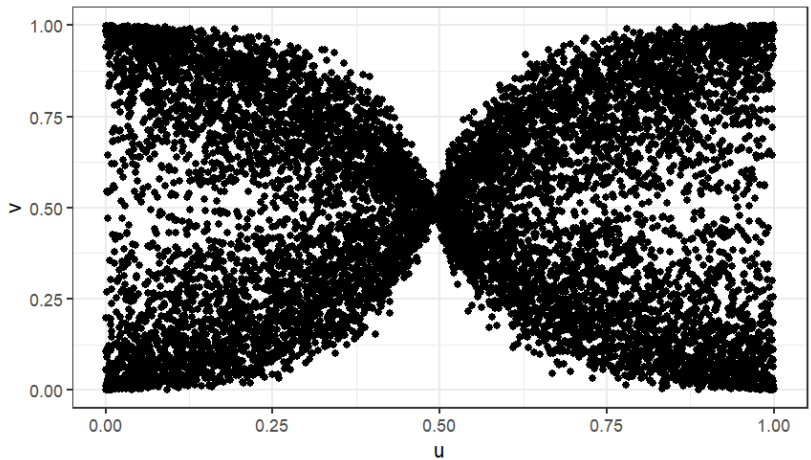

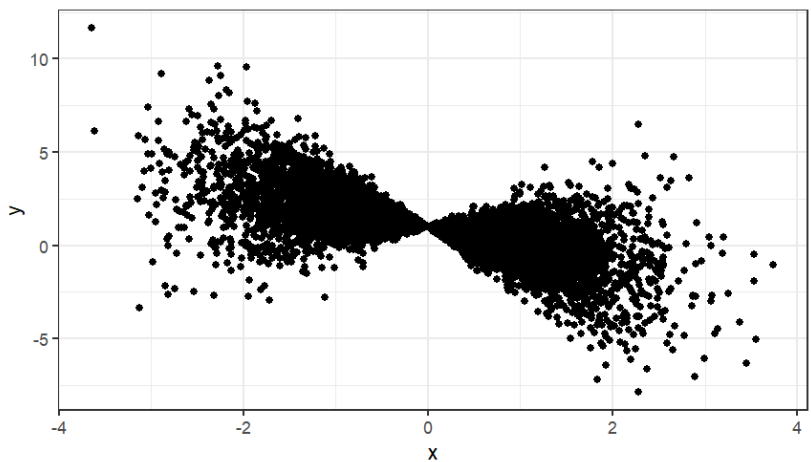

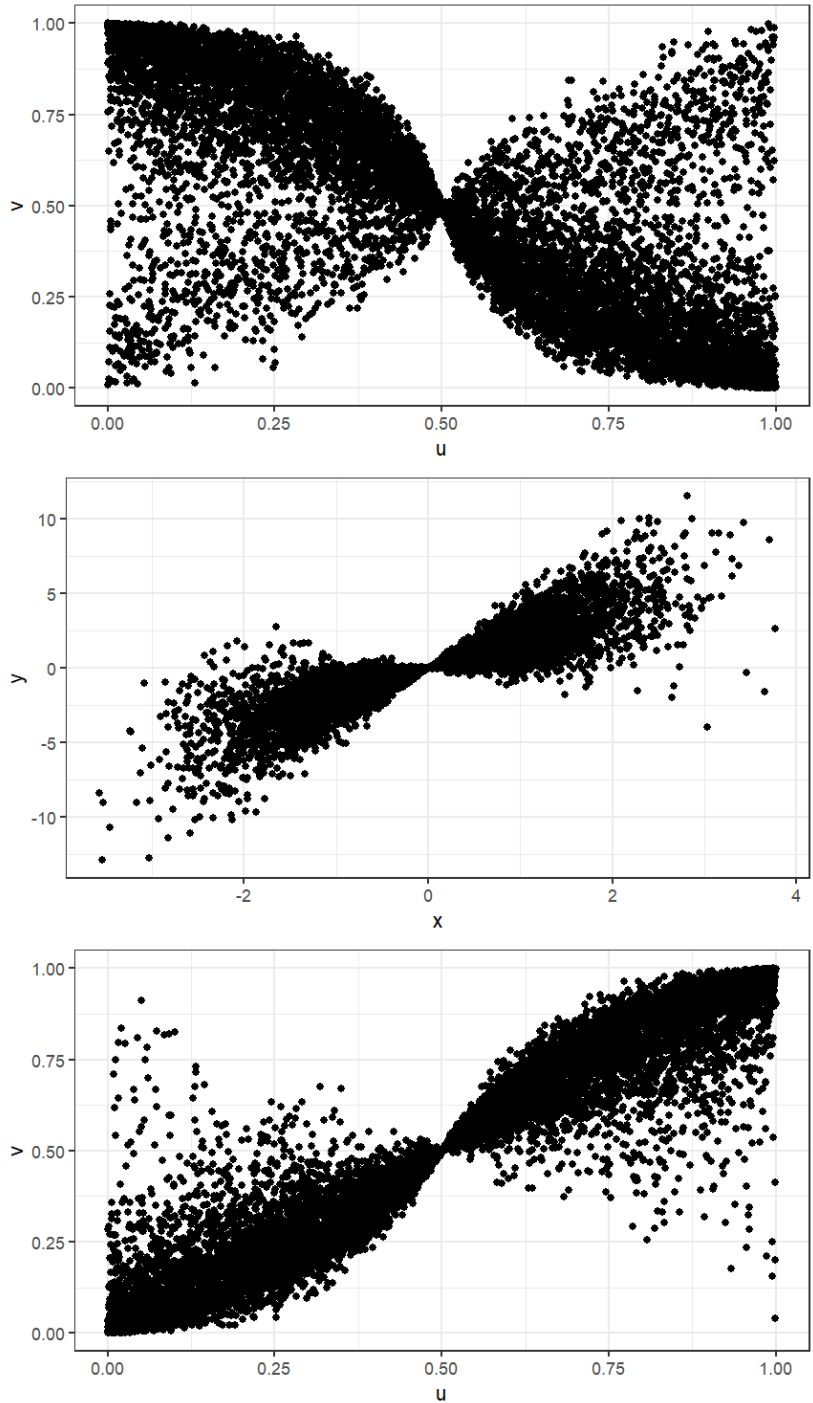

Figure 1: Scatter plots of $n = 10000$ pairs $(X, Y)$ generated from the cdf defined in Example 1 with the parameters $(\beta_0, \beta_1) = (0,0)$ (left panel top), $(\beta_0, \beta_1) = (1, -1)$ (left panel middel), $(\beta_0, \beta_1) = (0,2)$ (left panel bottom) and their corresponding normalized ranks (the right panels).

One approach to measuring dependence between two random variables in different regions of their joint distribution is to use the notion of quantile dependence. In this paper, we follow a different approach to measure the degree of dependence between quantiles, by using the conditional probability to develop an extension of tail dependence coefficients as a quantile dependence coefficient to detect the strength of dependence between two random variables in specific parts of their distribution. The paper is organized as follows. In Section 2, we present the main ideas related to the newly introduced index. The properties of the proposed index are studied in this section. Section 3 provided several examples. Section 4 concludes the paper.

## 2    The proposed quantile dependence coefficient and its properties

In this section, we proposed a family of quantile dependence coefficients and investigate their properties.

### 2.1    The proposed quantile dependence coefficient

Let $(X, Y)$ be a pair of continuous random variables with the joint distribution function $H$ and marginal distribution functions $F$ (of $X$) and $G$ (of $Y$). For $p, q \in [0,1]$, let $F^{-1}(p)$ and $G^{-1}(q)$ be the quantile function of $X$ and $Y$, respectively. Our proposed family of quantile dependence coefficients is given in Definition 2.

**Definition 2**    *For $p, q \in [0,1]$, we define the $(p, q)$-quantile dependence coefficient of a pair $(X, Y)$ by*

$$\lambda_{Y|X}(q|p) = \lim_{t \to 0^+} P\{G^{-1}((q-t)^+) < Y \leq G^{-1}((q+t)^-)|F^{-1}((p-t)^+) < X \leq F^{-1}((p+t)^-)\}, \tag{2.1}$$

*and*

$$\lambda_{X|Y}(p|q) = \lim_{t \to 0^+} P\{F^{-1}((p-t)^+) < X \leq F^{-1}((p+t)^-)|G^{-1}((q-t)^+) < Y \leq G^{-1}((q+t)^-)\}, \tag{2.2}$$

*if the limit exist, where $a^+ = max(a, 0)$ and $a^- = 1 - (1-a)^+$.*

The $(p, q)$-quantile dependence coefficient indicates the limit of the probability that one of the variables falls in an interval quantile of its distribution given that the other variable falls in the interval quantile of its distribution. By analyzing these limits in different regions of the joint distribution, one can gain insights into the degree of dependence between variables in those regions. We note that

$$\lambda_{Y|X}(0|0) = \lambda_{X|Y}(0|0))$$
$$= \lim_{t \to 0^+} P\{Y \leq G^{-1}(t)|X \leq F^{-1}(t)\} = \lambda_L,$$

and

$$\lambda_{Y|X}(1|1) = \lambda_{X|Y}(1|1)$$
$$= \lim_{t \to 0^+} P\{Y > G^{-1}(1-t)|X > F^{-1}(1-t)\}$$
$$= \lim_{t \to 1^-} P\{Y > G^{-1}(t)|X > F^{-1}(t)\} = \lambda_U.$$

The calculation of $\lambda_{Y|X}(q|p)$ and $\lambda_{X|Y}(p|q)$ can be simplified if the distribution of $(X, Y)$ is represented in terms of the copula. The following proposition defines the $(p, q)$-quantile

dependence coefficient via the notion of the copula.

**Remark 1** *For a pair $(X, Y)$ of continuous random variables, with the copula $C$, we use the notions $\lambda^C_{Y|X}(q|p)$ and $\lambda^C_{X|Y}(p|q)$, for $(p,q)$-quantile dependence coefficients of $(X, Y)$.*

**Proposition 1** *Let $(X, Y)$ be a pair of continuous random variables with the joint distribution function $H$ and marginal distribution functions $F$ (of $X$) and $G$ (of $Y$). If $C$ is the copula of $(X, Y)$ then for $p, q \in [0,1]$*

$$\lambda^C_{Y|X}(q|p) = \lim_{t \to 0^+} \frac{V_C([(p-t)^+, (p+t)^-] \times [(q-t)^+, (q+t)^-])}{(p+t)^- - (p-t)^+},$$

*and*

$$\lambda^C_{X|Y}(p|q) = \lim_{t \to 0^+} \frac{V_C([(p-t)^+, (p+t)^-] \times [(q-t)^+, (q+t)^-])}{(q+t)^- - (q-t)^+},$$

*provided that the limits exist, where for $0 \leq u_1 \leq u_2 \leq 1$ and $0 \leq v_1 \leq v_2 \leq 1$, $V_C([u_1, u_2] \times [v_1, v_2]) = C(u_2, v_2) - C(u_2, v_1) - C(u_1, v_2) + C(u_1, v_1)$ is the $C$-volume of the rectangle $[u_1, u_2] \times [v_1, v_2]$.*

*Proof.* For $p, q \in [0,1]$ we have that

$$P[F^{-1}((p-t)^+) < X \leq F^{-1}((p+t)^-)] = F(F^{-1}(p+t)) - F(F^{-1}(p-t))$$
$$= (p+t)^- - (p-t)^+.$$

Therefore

$$\lambda^C_{Y|X}(q|p) = \lim_{t \to 0^+} P[G^{-1}((q-t)^+) < Y \leq G^{-1}((q+t)^-)|F^{-1}((p-t)^+) < X$$
$$\leq F^{-1}((p+t)^-)]$$
$$= \lim_{t \to 0^+} \frac{P[F^{-1}((p-t)^+) < X \leq F^{-1}((p+t)^-), G^{-1}((q-t)^+) < Y \leq G^{-1}((q+t)^-)]}{P[F^{-1}(p-t)^+ < X \leq F^{-1}(p+t)^-]}$$
$$= \lim_{t \to 0^+} \frac{1}{(p+t)^- - (p-t)^+} [H(F^{-1}((p+t)^-), G^{-1}((q+t)^-))$$
$$- H(F^{-1}((p+t)^-), G^{-1}((q-t)^+)) - H(F^{-1}((p-t)^+), G^{-1}((q+t)^-))$$
$$+ H(F^{-1}((p-t)^+), G^{-1}((q-t)^-))]$$
$$= \lim_{t \to 0^+} \frac{1}{(p+t)^- - (p-t)^+} [C((p+t)^-, (q+t)^-) - C((p+t)^-, (q-t)^+)$$
$$- C((p-t)^+, (q+t)^-) + C((p-t)^+, (q-t)^+)],$$

where the later result follows from Sklar's Theorem that $H(F^{-1}(u), G^{-1}(v)) = C(u, v)$. A similar argument holds for $\lambda^C_{X|Y}(p|q)$. □

We note that

$$\lambda^C_{Y|X}(0|0) = \lambda^C_{X|Y}(0|0) = \lim_{t \to 0^+} \frac{C(t,t)}{t} = \lambda_L.$$

and

$$\lambda^C_{Y|X}(1|1) = \lambda^C_{X|Y}(1|1) = \lim_{t \to 0^+} \frac{2t - 1 + C(1-t, 1-t)}{t} = \lambda_U.$$

Thus $\lambda^C_{X|Y}(p|q)$ and $\lambda^C_{Y|X}(q|p)$ are extensions of tail dependence coefficients $\lambda_L$ and $\lambda_U$. From the Definition 2, for every copula $C$ and $p, q \in [0,1]$, it follows that $\lambda^C_{Y|X}(q|p) \in [0,1]$ and $\lambda^C_{X|Y}(p|q) \in [0,1]$. For the copula $\Pi(u, v) = uv$ of independent random variables, it is easy to see that $\lambda^\Pi_{X|Y}(p|q) = \lambda^\Pi_{Y|X}(q|p) = 0$, for all $p, q \in [0,1]$. If $\lambda^C_{Y|X}(q|p) \in (0,1]$ or $\lambda^C_{X|Y}(p|q) \in (0,1]$, for some $p, q \in [0,1]$, we say $C$ has $(p,q)$-quantile dependence; if $\lambda^C_{Y|X}(q|p) = 0$ or $\lambda^C_{X|Y}(p|q) = 0$, we say $C$ has no $(p,q)$-quantile dependence.

The following Corollary provides the expression of $\lambda^C_{Y|X}(q|p)$ for different values of $p$ and $q$.

**Corollary 1**  *For $p, q \in (0,1)$*

$$\lambda^C_{Y|X}(q|p) = \lim_{t \to 0^+} \frac{C(p+t,q+t) - C(p+t,q-t) - C(p-t,q+t) + C(p-t,q-t)}{2t},$$

*and*

$$\lambda^C_{Y|X}(q|p) = \lim_{t \to 0^+} \begin{cases} \frac{C(t,q+t) - C(t,q-t)}{t}, & p = 0, q \in (0,1) \\ \frac{2t - C(1-t,q+t) + C(1-t,q-t)}{t}, & p = 1, q \in (0,1), \\ \frac{C(p+t,t) - C(p-t,t))}{2t}, & q = 0, p \in (0,1), \\ \frac{2t - C(p+t,1-t) - C(p-t,1-t))}{2t}, & q = 1, p \in (0,1), \\ \frac{t - C(1-t,t)}{t}, & q = 0, p = 1, \\ \frac{t - C(t,1-t)}{t}, & q = 1, p = 0. \end{cases}$$

**Remark 2** *Note that in general $\lambda^C_{Y|X}(q|p)$ and $\lambda^C_{X|Y}(p|q)$ are not equal, but the following relations hold between them:*

$$\lambda^C_{X|Y}(p|q) = \begin{cases} \lambda^C_{Y|X}(q|p), & p, q \in (0,1) \quad or \quad p, q \in \{0,1\}, \\ 2\lambda^C_{Y|X}(q|p), & p \in (0,1), q \in \{0,1\}, \\ \frac{1}{2}\lambda^C_{Y|X}(q|p), & p \in \{0,1\}, q \in (0,1). \end{cases}$$

*Therefore,*

$$\lambda^C_{Y|X}(q|p) \leq \frac{1}{2}, \qquad for \qquad q \in \{0,1\}, p \in (0,1),$$

*and*

$$\lambda^C_{X|Y}(p|q) \leq \frac{1}{2}, \qquad for \qquad q \in (0,1), p \in \{0,1\}.$$

In the following example, we provide the values of the $(p,q)$-quantile measure of dependence for the regression model given in Example 1.

Table   1: The values of $\lambda_{Y|X}(\frac{1}{2}|\frac{1}{2})$ for the model given in Example 1

| $\beta_1$ | 0 | $\pm 1$ | $\pm 2$ | $\pm 3$ |
|---|---|---|---|---|
| $\lambda_{Y|X}(\frac{1}{2}, \frac{1}{2})$ | 0.1823 | 0.1795 | 0.3122 | 0.5 |

**Example 2**  *Consider the regression model given in Example 1. The joint cdf of $(X,Y)$ is given by*

$$H(x,y) = \int_{-\infty}^{x} \int_{-\infty}^{y} f_X(s) f_{Y|X=s}(t) dt ds$$
$$= \int_{-\infty}^{x} \phi(s) \Phi(\frac{y - \beta_0 - \beta_1 s}{|s|}) ds,$$

where $\Phi(.)$ and $\phi(.)$ are the cumulative distribution function and the density function of $N(0,1)$ random variable. In view of Sklar's Theorem the copula of $(X, Y)$ is given by

$$C(u, v) = \int_0^u \Phi^{-1}\left(\frac{G^{-1}(v) - \beta_1 \Phi^{-1}(w)}{|\Phi^{-1}(w)|}\right) dw, \qquad (2.3)$$

where $G^{-1}(.)$ is the inverse of the cdf of $Y - \beta_0$, given by

$$G_Y(y) = \int_{-\infty}^{\infty} \phi(x) \Phi\left(\frac{y - \beta_1 x}{|x|}\right) dx.$$

Note that the marginal cdf of $Y$ is symmetric about $\beta_0$. A straightforward calculation shows that $(X, Y)$ has $(\frac{1}{2}, \frac{1}{2})$-quantile or median dependence, as we see in Figure 1. We note that the value of $(\frac{1}{2}, \frac{1}{2})$ does not depend on the parameter $\beta_0$ and its value is the same for $-\beta_1$ and $+\beta$. Table 1 shows the values of $\lambda_{Y|X}(\frac{1}{2}|\frac{1}{2})$ for different values of the parameter $\pm\beta_1$. As we saw in Example 1, for the case $\beta_1 = 0$ the value of the Pearson's correlation coefficient of this model is equal to zero but $\lambda_{Y|X}(\frac{1}{2}|\frac{1}{2}) = 0.1823$.

The following example shows a copula for which, in addition to the upper and lower tail dependency, there is dependency in other regions of the support.

**Example 3** Let $X \sim U(0,1)$ and $Y = 3X\ I_{(0,\frac{1}{3})} - (3X - 2)\ I_{[\frac{1}{3},\frac{2}{3}]} + (3X - 2)\ I_{[\frac{2}{3},1)}$. By straightforward calculation the copula of $(X, Y)$ is given by

$$C(u, v) = \begin{cases} u, & 0 < u < \frac{v}{3} \\ \frac{v}{3}, & \frac{v}{3} \le u < \frac{-v+2}{3} \\ u - \frac{2}{3}(1 - v), & \frac{-v+2}{3} \le u < \frac{v+2}{3} \\ v, & \frac{v+2}{3} \le u < 1 \end{cases} \qquad (2.4)$$

The quantile dependence coefficient $\lambda_{Y|X}^C(q|p)$ is then

$$\lambda_{Y|X}^C(q|p) = \begin{cases} \frac{1}{3}, & p = \frac{q}{3}, p = \frac{-q+2}{3},\ p = \frac{q+2}{3} \\ 0, & o.w \end{cases} \qquad (2.5)$$

Figure 2 shows the plot of $\lambda_{Y|X}^C(q|p)$ for this copula.

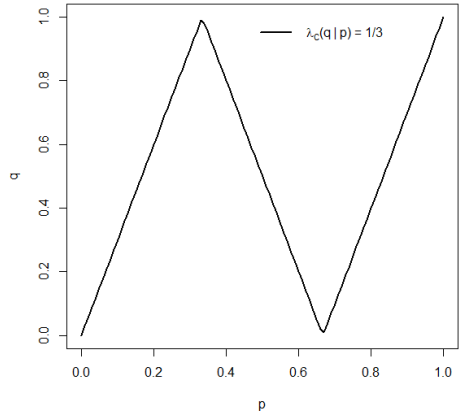

Figure 2: Plot of $(p, q)$-quantile dependence coefficient of the copula defined in Example 3

## 2.2    Properties of the proposed quantile dependence coefficient

In this section, we discuss the properties of the proposed index. First we consider $(p, q)$-quantile dependence coefficient for the Fréchet-Hoeffding upper bound copula $M(u, v) = \min(u, v)$ and the Fréchet-Hoeffding lower bound copula $W(u, v) = \max(u + v - 1, 0)$ [11].

**Remark 3** *In the following, for the sake of simplicity, we use the notions $\lambda^C(q|p)$ and $\lambda^C(p|q)$ instead of $\lambda^C_{Y|X}(q|p)$ and $\lambda^C_{X|Y}(p|q)$, respectively. We will study properties of $\lambda^C(q|p)$, but similar results hold for $\lambda^C(p|q)$.*

**Proposition 2**    *For the copulas $M$ and $W$ we have*
$$\lambda^M(q|p) = I_{\{p=q\}}, \qquad and \qquad \lambda^W(q|p) = I_{\{p+q=1\}},$$
*for all $p, q \in [0,1]$.*
*Proof.* The proof is straightforward from Proposition 1.                                      □

Observe that the convex combination of two copulas is a copula [11]. The following result shows that the quantile dependence coefficient for a convex combination of two copulas is the convex combination of their quantile dependence coefficients.

**Proposition 3** *For $\omega \in [0,1]$, let $C(u, v) = \omega C_1(u, v) + (1 - \omega)C_2(u, v)$, where $C_1$ and $C_2$ are two copulas. Then for all $p, q \in [0,1]$,*
$$\lambda^C(q|p) = \omega \lambda^{C_1}(q|p) + (1 - \omega)\lambda^{C_2}(q|p).$$
*Proof.* The proof is straightforward from Proposition 1 using the fact that the $C$-volume of the rectangle $B \subset [0,1] \times [0,1]$ is given by $V_C(B) = \omega V_{C_1}(B) + (1 - \omega)V_{C_2}(B)$.                  □

**Example 4**    *For $\omega \in [0,1]$, let $A(u, v) = \omega M(u, v) + (1 - \omega)W(u, v)$. Then*
$$\lambda^A(q|p) = \omega \ I_{(p=q)} + (1 - \omega) \ I_{(p=1-q)}.$$
*Figure 5 shows the plot of $(p, q)$-quantile dependence coefficient of the copula A for $\omega = \frac{2}{3}$.*

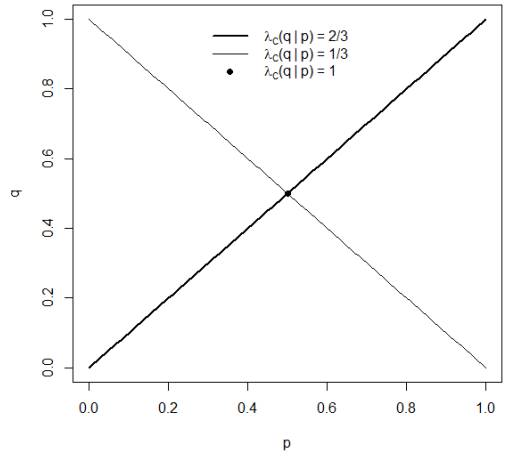

Figure 3: Plot of $(p, q)$-quantile dependence coefficient of the copula A defined in Example 4 for $\omega = \frac{2}{3}$.

For the copula $A$ in Example 4, $\lambda^A(\frac{1}{2}|\frac{1}{2}) = 1$ and $\lambda^A(q|\frac{1}{2}) = 0$, for every $q \in [0,1] - \{\frac{1}{2}\}$. Therefore, $\lambda^A(\frac{1}{2}|\frac{1}{2}) + \lambda^A(q|\frac{1}{2}) = 1$, for all $q \neq \frac{1}{2}$. The following result shows a similar result for every copula $C$.

**Proposition 4** *Let* $p_0, q_0$ *be two points in* $[0,1]$. *Then, for every copula* $C$ *and* $q \in [0,1]$ *with* $q \neq q_0$,
$$\lambda^C(q_0|p_0) + \lambda^C(q|p_0) \leq 1.$$
*Proof.* Let $(U,V) \sim C$ and $A_1 = \{V \in ((q_0-t)^+, (q_0+t)^-]\}$, $A_2 = \{V \in ((q-t)^+, (q+t)^-]\}$ and $B = \{U \in ((p_0-t)^+, (p_0+t)^-]\}$. Then from Definition 2,
$$\lambda^C(q_0|p_0) + \lambda^C(q|p_0) = \lim_{t \to 0^+} \frac{P(A_1 \cap B) + P(A_2 \cap B)}{P(B)}.$$
For $t < \frac{|q_0-q|}{2}$, we have $A_1 \cap A_2 = \emptyset$ and thus $P(A_1 \cap B) + P(A_2 \cap B) = P(A \cap B) \leq P(B)$, where $A = A_1 \cup A_2$, which completes the proof. $\qquad \square$

Note that we can also derive tail dependence coefficients of a survival copula from its associated copula using the (1.3) and (1.4) via $\lambda_L(\hat{C}) = \lambda_U(C)$ and $\lambda_U(\hat{C}) = \lambda_L(C)$; that is, for a copula with an upper tail dependence, its associated survival copula has a lower tail dependence and conversely. In the following, we provide some results for symmetry properties of $(p,q)$-quantile dependence coefficient.

**Proposition 5** *Let* $\hat{C}$ *be the survival copula associated with the copula* $C$. *Then,*
$$\lambda^{\hat{C}}(q|p) = \lambda^C(1-q|1-p).$$
*Proof.* The proof is straightforward from Proposition 1. $\qquad \square$

We recall that if $C$ is the copula of the pair $(U,V)$, then the reflections of the copula $C$ given by
$$C^*(u,v) = u - C(u, 1-v), \qquad \text{and} \qquad C^{**}(u,v) = v - C(1-u, v), \qquad (2.6)$$
are the copulas of the pairs $(U, 1-V)$ and $(1-U, V)$, respectively (see, e.g, [11], Theorem 2.4.4). A copula $C$ is called conditionally symmetric if $C(u,v) = C^*(u,v)$ or $C(u,v) = C^{**}(u,v)$, for all $u,v \in [0,1]$. If $C^*(u,v) = C^{**}(u,v) = C(u,v)$, for all $u,v \in [0,1]$, then $C$ is called jointly symmetric [11]. Jointly symmetric random variables must be uncorrelated when their second-order moments exist. In fact, for jointly symmetric random variables, all of the concordance measures satisfy Scarsini's axioms [14], such as Kendall's tau, Spearman's rho, and Gini's gamma are equal to zero [2]. The following result shows that we can derive the quantile measure of dependence of the copulas $C^*$ and $C^{**}$, from their associated copula $C$.

**Proposition 6** *For a given copula* $C$, *let* $C^*$ *and* $C^{**}$ *be the associated copulas given by (2.6). Then,*
$$\lambda^{C^*}(q|p) = \lambda^C(1-q|p), \qquad \text{and} \qquad \lambda^{C^{**}}(q|p) = \lambda^C(q|1-p). \qquad (2.7)$$
*Proof.* The result follows from Corollary 1. $\qquad \square$

**Remark 4** *We note that for a jointly symmetric copula* $C$, *that is* $C^* = C^{**} = C$, *we have that for all* $p, q \in [0,1]$,

$$\lambda^C(1-q|p) = \lambda^C(q|1-p) = \lambda^C(1-q|1-p) = \lambda^C(q|p).$$

**Example 5** *Consider the copula* $A(u,v) = \frac{M(u,v)+W(u,v)}{2}$. *This copula is jointly symmetric and from Example 4,* $\lambda^A(1-p|p) = \lambda^A(p|1-p) = \lambda^A(1-p|1-p) = \lambda^A(p|p) = \frac{1}{2}$ *for all* $p \in [0,1] - \{\frac{1}{2}\}$.

The following proposition shows the continuity property of the proposed quantile dependence coefficient. That is, if the copula of a sequence of bivariate random pairs converges in distribution to $C$, then the sequence of the coefficients for these random pairs converges to the measure of $C$.

**Proposition 7** *Let* $\{C_n\}_{n \in \, IN}$ *be a sequence of copulas such that* $C_n \to C$, *then*
$$lim_{n \to \infty} \lambda^{C_n}(q|p) = \lambda^C(q|p).$$

*Proof.* Since $C_n \to C$, from the Proposition 1 we have that
$$\lim_{n \to \infty} \lambda^{C_n}(q|p) = \lim_{n \to \infty} \lim_{t \to 0+} \frac{V_{C_n}([(p-t)^+,(p+t)^-] \times [(q-t)^+,(q+t)^-])}{(p+t)^- - (p-t)^+}$$
$$= \lim_{t \to 0+} \lim_{n \to \infty} \frac{V_{C_n}([(p-t)^+,(p+t)^-] \times [(q-t)^+,(q+t)^-])}{(p+t)^- - (p-t)^+}$$
$$= \lim_{t \to 0+} \frac{V_C([(p-t)^+,(p+t)^-] \times [(q-t)^+,(q+t)^-])}{(p+t)^- - (p-t)^+} = \lambda^C(q|p). \qquad \square$$

## 3   Quantile dependence coefficient for some families of copulas

In this section, we examine the value of the proposed $(p,q)$-quantile dependence coefficient for some families of copulas such as Normal copula, Extreme value copula, and Archimedean copulas. First, we show that the proposed quantile dependence coefficient $\lambda^C(.|.)$ can also be calculated using the conditional cdfs of the copula $C$. For a pair $(U,V)$ of uniform $(0,1)$ random variables with the copula $C$, the partial derivative $C_{2|1}(v|u) = \partial C(u,v)/\partial u$ is the conditional cdf of $[V|U = u]$ and the partial derivative $C_{1|2}(u|v) = \partial C(u,v)/\partial v$ is the conditional cdf of $[U|V = v]$. For any $v \in [0,1]$, the partial derivative $C_{2|1}(v|u)$ exists for almost all $u$, and for such $v$ and $u$, $C_{2|1}(v|u) \in [0,1]$. Similarly, for any $u \in [0,1]$, the partial derivative $C_{1|2}(u|v)$ exists for almost all $\in [0,1]$, and for such $u$ and $v$, $C_{1|2}(u|v) \in [0,1]$; see, e.g., [11]. The following useful result shows that $(p,q)$-quantile dependence coefficient can be studied through the conditional cdfs $C_{2|1}(v|u)$ and $C_{1|2}(u|v)$.

**Proposition 8**   *Let* $\lambda^C(q|p)$ *be as in Proposition 1. If the limits exist, then*
$$\lambda^C(q|p) = \frac{\zeta_{2|1}(q|p) + \zeta_{1|2}(q|p)}{1 + I_{[p \in (0,1)]}},$$

*where*
$$\zeta_{2|1}(q|p) = \lim_{t \to 0+} I_{[p+t<1]}[C_{2|1}((q+t)^-|(p+t)^-) - C_{2|1}((q-t)^+|(p+t)^-)]$$
$$+ \lim_{t \to 0+} I_{[p-t>0]}[C_{2|1}((q+t)^-|(p-t)^+) - I_{[p-t>0]}C_{2|1}((q-t)^+|(p-t)^+)],$$

*and*

$$\zeta_{1|2}(q|p) = \lim_{t\to 0^+} [\ I_{[q+t<1]}[C_{1|2}((p+t)^-|(q+t)^-) - I_{[q-t>0]}C_{1|2}((p+t)^-|(q-t)^+)]$$
$$+ \lim_{t\to 0^+}\ I_{[q+t<1]}[C_{1|2}((p-t)^+|(q+t)^-) - I_{[q-t>0]}C_{1|2}((p-t)^+|(q-t)^+)].$$

*Proof.* We note that

$$\lim_{t\to 0^+} \frac{d}{dt}[(p+t)^- - (p-t)^+] = \lim_{t\to 0^+} [\ I_{[p+t<1]} + \ I_{[p-t>0]}] = 1 + I_{[p\in(0,1)]}.$$

Since

$$\frac{d}{dt} C(g(t), h(t)) = g'(t)\frac{\partial C(u,v)}{\partial u}\Big|_{u=g(t),v=h(t)} + h'(t)\frac{\partial C(u,v)}{\partial v}\Big|_{u=g(t),v=h(t)},$$

then,

$$\frac{\partial}{\partial t} C((p+t)^-, (q+t)^-) = I_{[p+t<1]}C_{2|1}((q+t)^-|(p+t)^-) + I_{[q+t<1]}C_{1|2}((p+t)^-|(q+t)^-),$$

$$\frac{\partial}{\partial t} C((p+t)^-, (q-t)^+) = I_{[p+t<1]}C_{2|1}((q-t)^+|(p+t)^-) - I_{[q-t>0]}C_{1|2}((p+t)^-|(q-t)^+),$$

$$\frac{\partial}{\partial t} C((p-t)^+, (q+t)^-) = -I_{[p-t>0]}C_{2|1}((q+t)^-|(p-t)^+) + I_{[q+t<1]}C_{1|2}((p-t)^+|(q+t)^-),$$

and

$$\frac{\partial}{\partial t} C((p-t)^+, (q-t)^+) = -\ I_{[p-t>0]}C_{2|1}((q-t)^+|(p-t)^+) - I_{[q-t>0]}C_{1|2}((p-t)^+|(q-t)^+).$$

By using the technique of l'Hopital's we obtain the result. $\qquad\square$

Here are some special cases of $\lambda^C(q|p)$ in Proposition 8.

**Corollary 2** *For $p, q \in (0,1)$ we have that*
$$\lambda^C(q|p) = \frac{1}{2}\lim_{t\to 0^+}[C_{2|1}(q+t|p+t) - C_{2|1}(q-t|p+t) + C_{2|1}(q+t|p-t)$$
$$-C_{2|1}(q-t|p-t) + C_{1|2}(p+t|q+t) - C_{1|2}(p-t|q+t)$$
$$+C_{1|2}(p+t|q-t) - C_{1|2}(p-t|q-t)],$$

*for $p \in (0,1)$ and $q = 0$,*
$$\lambda^C(0|p) = \frac{1}{2}\lim_{t\to 0^+}[C_{2|1}(t|p+t) + C_{2|1}(t|p-t) + C_{1|2}(p+t|t) - C_{1|2}(p-t|t)],$$

*for $p \in (0,1)$ and $q = 1$.*
$$\lambda^C(1|p) = 1 - \frac{1}{2}\lim_{t\to 0^+}[C_{2|1}(1-t|p+t) + C_{2|1}(1-t|p-t) - C_{1|2}(p+t|1-t) + C_{1|2}(p-t|1-t)],$$

*for $p = 0$ and $q \in (0,1)$,*
$$\lambda^C(q|0) = \lim_{t\to 0^+}[C_{2|1}(q+t|t) - C_{2|1}(q-t|t) + C_{1|2}(t|q+t) + C_{1|2}(t|q-t)],$$

*for $p = 1$ and $q \in (0,1)$,*
$$\lambda^C(q|0) = 2 - \lim_{t\to 0^+}[C_{1|2}(1-t|q+t) + C_{1|2}(1-t|q-t) - C_{2|1}(q+t|1-t) + C_{2|1}(q-t|1-t)],$$

*for $p = 0$ and $q = 1$,*
$$\lambda^C(1|0) = 1 - \lim_{t\to 0^+}[C_{2|1}(1-t|t) - C_{1|2}(t|1-t)],$$

*and for $q = 0$ and $p = 1$,*
$$\lambda^C(0|1) = 1 - \lim_{t\to 0^+}[C_{1|2}(1-t|t) - C_{2|1}(t|1-t)].$$

Note that
$$\lambda_L = \lambda^C(0|0) = \lim_{t\to 0^+}[C_{2|1}(t|t) + C_{1|2}(t|t)],$$

and
$$\lambda_U = \lambda^C(1|1) = 2 - \lim_{t \to 0^+}[C_{2|1}(1-t|1-t) + C_{1|2}(1-t|1-t)].$$

## 3.1 Gaussian Copula

In this subsection, we study the quantile dependence coefficient of the bivariate Gaussian (Normal) copula. The Gaussian copula $C^\rho(u,v)$ is defined, for $u,v \in [0,1]$ as follows:
$$C^\rho(u,v) = \Phi_\rho(\Phi^{-1}(u), \Phi^{-1}(v)),$$

where $\Phi$ and $\Phi^{-1}$ are the standard-normal distribution function and its inverse, and $\Phi_\rho$ is the cdf of the bivariate normal distribution with the correlation parameter $-1 \le \rho \le 1$. For this copula
$$C^\rho_{2|1}(v|u) = \Phi\left(\frac{\Phi^{-1}(v) - \rho\Phi^{-1}(u)}{\sqrt{1-\rho^2}}\right),$$

and
$$C^\rho_{1|2}(u|v) = \Phi\left(\frac{\Phi^{-1}(u) - \rho\Phi^{-1}(v)}{\sqrt{1-\rho^2}}\right).$$

The functions $\Phi(t)$ and $\Phi^{-1}(t)$ are continuous on IR and $(0,1)$, respectively. Thus, the functions $C^\rho_{1|2}(u|v)$ and $C^\rho_{2|1}(v|u)$ are continuous in $u$ and $v$ on $(0,1)$. These functions are also right-continuous in $u$ and $v$ at zero and left-continuous at one. The lower and upper tail dependence coefficients for Gaussian copula is given by $\lambda_L = \lambda_U = 0$. For every $\rho \in (-1,1)$ and every $p,q \in (0,1)$, from Corollary 2, we have that
$$\lim_{t \to 0^+} C^\rho_{2|1}(q+t|p+t) = \lim_{t \to 0^+} C^\rho_{2|1}(q-t|p+t) = \lim_{t \to 0^+} C^\rho_{2|1}(q+t|p-t)$$
$$= \lim_{t \to 0^+} C^\rho_{2|1}(q-t|p-t) = C^\rho_{2|1}(q|p),$$

and
$$\lim_{t \to 0^+} C^\rho_{1|2}(p+t|q+t) = \lim_{t \to 0^+} C^\rho_{1|2}(p-t|q+t) = \lim_{t \to 0^+} C^\rho_{1|2}(p+t|q-t)$$
$$= C^\rho_{1|2}(p-t|q-t) = C^\rho_{1|2}(p|q),$$

and thus for all $\rho \in (-1,1)$, $\lambda^{C^\rho}(q|p) = 0$ for all $p,q \in [0,1]$. Therefore, for every $\rho \in (-1,1)$, the Gaussian copula is $(p,q)$-quantile independence for all $p,q \in [0,1]$. Since, $\lim_{\rho \to +1} C^\rho(u,v) = M(u,v)$ and $\lim_{\rho \to -1} C^\rho(u,v) = W(u,v)$, then from Proposition 2, we have that $\lambda^{C^{+1}}(q|p) = I_{[p=q]}$ and $\lambda^{C^{-1}}(q|p) = I_{[p+q=1]}$. Thus, the Gaussian copula has $(p,q)$-quantile dependence for all $p,q \in [0,1]$, only for the case $\rho \in \{-1,1\}$.

## 3.2 Student T copula

In this subsection, we study the quantile dependence coefficient of the bivariate Student T copula. The Student T copula $C^{\rho,\nu}(u,v)$ is defined, for $u,v \in [0,1]$ as follwos:
$$C^{\rho,\nu}(u,v) = T_{\rho,\nu}(T_\nu^{-1}(u), T_\nu^{-1}(v)),$$

where $T_\nu(.)$ and $T_\nu^{-1}(.)$ are the cdf of Student T random variable with the degrees of freedom $\nu \in \mathbb{N} - \{0\}$ and its inverse, and $T_{\rho,\nu}(.,.)$ is the cdf of the bivariate Student T distribution with the correlation parameter $-1 \le \rho \le 1$ and the degrees of freedom $\nu$. Note that $\lim_{\rho \to +1} C^{\rho,\nu}(u,v) = M(u,v)$ and $\lim_{\rho \to -1} C^{\rho,\nu}(u,v) = W(u,v)$. For this copula

$$C_{2|1}^{\rho,\nu}(v|u) = T_{\nu+1}\left(\frac{T_\nu^{-1}(v) - \rho T_\nu^{-1}(u)}{\sqrt{(1-\rho^2)(\nu + [T_\nu^{-1}(u)]^2)/(\nu+1)}}\right),$$

and

$$C_{1|2}^{\rho,\nu}(u|v) = T_{\nu+1}\left(\frac{T_\nu^{-1}(u) - \rho T_\nu^{-1}(v)}{\sqrt{(1-\rho^2)(\nu + [T_\nu^{-1}(v)]^2)/(\nu+1)}}\right).$$

The functions $T_\nu(.)$ and $T_\nu^{-1}(.)$ are continuous on $\mathbb{IR}$ and $(0,1)$, respectively. Thus, the functions $C_{1|2}^{\rho,\nu}(u|v)$ and $C_{2|1}^{\rho,\nu}(v|u)$ are continuous in $u$ and $v$ on $(0,1)$. From Corollary 2 the $(p,q)$-quantile dependence coefficient $\lambda^{C^{\rho,\nu}(q|p)}$ for the Student T copula is given as follows:

$$\lambda^{C^{\rho,\nu}}(0|0) = \lambda^{C^{\rho,\nu}}(1|1) = 2T_{\nu+1}\left(-\sqrt{\frac{(\nu+1)(1-\rho)}{1+\rho}}\right),$$

$$\lambda^{C^{\rho,\nu}}(1|0) = \lambda^{C^{\rho,\nu}}(0|1) = 2T_{\nu+1}\left(-\sqrt{\frac{(\nu+1)(1+\rho)}{1-\rho}}\right).$$

For every $\rho \in (-1,1)$ and every $p,q \in (0,1)$, from Corollary 2, it is easy to see that $\lambda^{C^{\rho,\nu}}(q|p) = \lambda^{C^{\rho,\nu}}(q|1) = \lambda^{C^{\rho,\nu}}(q|0) = \lambda^{C^{\rho,\nu}}(0|p) = \lambda^{C^{\rho,\nu}}(1|p) = 0.$

### 3.3   Extreme Value copulas

Let $C$ be an extreme value copula (EV) defined by [11]

$$C(u,v) = \exp\left(\ln(uv)A\left(\frac{\ln(v)}{\ln(uv)}\right)\right), \tag{3.1}$$

where $A:[0,1] \to [1/2,1]$ satisfies $A(0) = A(1) = 1$ and $\max(t, 1-t) \le A(t) \le 1$ is its dependence function. For an EV copula, the upper tail dependence coefficient is given by $\lambda_U = 2(1 - A(1/2))$, which can be interpreted as the length between the upper boundary and the curve $A(.)$ evaluated in the mid-point $1/2$. The coefficient $\lambda_U$ ranges from 0 (independence) to 1 (complete dependence). The lower tail dependence coefficient is given by $\lambda_L = I_{[A(1/2)=1/2]}$. That is, except for the case of perfect dependence, $A(1/2) = 1/2$, EV copulas have asymptotically independent lower tails [11].

The following proposition provides the $(p,q)$-quantile dependence coefficient of EV copulas.

**Proposition 9**   *Let $C$ be an extreme value copula given by (3.1). If the function $A(t)$ is continuous at the point $\Delta = \frac{\ln(q)}{\ln(pq)}$, for $p,q \in (0,1)$, then*

$$\lambda^C(q|p) = \frac{\ln(q)}{p\ln(pq)}\exp\left[\ln(pq)A\left(\frac{\ln(q)}{\ln(pq)}\right)\right]\left[\lim_{t\to\Delta^+}A'(t) - \lim_{t\to\Delta^-}A'(t)\right]. \tag{3.2}$$

*Proof.* The conditional cdfs $C_{2|1}(v|u)$ and $C_{1|2}(u|v)$ for $C$ are given by

$$C_{2|1}(v|u) = \frac{C(u,v)}{u}\left(A\left(\frac{\ln(v)}{\ln(uv)}\right) - \frac{\ln(v)}{\ln(uv)}A'\left(\frac{\ln(v)}{\ln(uv)}\right)\right),$$

and

$$C_{1|2}(u|v) = \frac{C(u,v)}{v}\left(A\left(\frac{\ln(v)}{\ln(uv)}\right) + \frac{\ln(u)}{\ln(uv)}A'\left(\frac{\ln(v)}{\ln(uv)}\right)\right).$$

Thus,

$$\lim_{t \to 0^+} C_{2|1}(q+t|p+t) = \lim_{t \to 0^+} C_{2|1}(q+t|p-t) = \lim_{t \to \Delta^-} Q_1(t),$$

$$\lim_{t \to 0^+} C_{2|1}(q-t|p+t) = \lim_{t \to 0^+} C_{2|1}(q-t|p-t) = \lim_{t \to \Delta^+} Q_1(t),$$

$$\lim_{t \to 0^-} C_{1|2}(p+t|q+t) = \lim_{t \to 0^+} C_{1|2}(p-t|q+t) = \lim_{t \to \Delta^-} Q_2(t),$$

$$\lim_{t \to 0^+} C_{1|2}(p+t|q-t) = \lim_{t \to 0^+} C_{1|2}(p-t|q-t) = \lim_{t \to \Delta^+} Q_2(t),$$

where,

$$Q_1(t) = \frac{1}{p} \exp\{A(t)\ln(pq)\}[A(t) - \frac{\ln(q)}{\ln(pq)} A'(t)],$$

and

$$Q_2(t) = \frac{1}{q} \exp\{A(t)\ln(pq)\}[A(t) - \frac{\ln(p)}{\ln(pq)} A'(t)].$$

Now from Proposition 8, if $A(t)$ is continuous at $t = \Delta(= \frac{\ln(q)}{\ln(pq)})$, then we have that

$$\lambda^C(q|p) = \lim_{t \to \Delta^-} Q_1(t) - \lim_{t \to \Delta^+} Q_1(t)$$
$$= \frac{\ln(q)}{p\ln(pq)} \exp\left[\ln(pq)A\left(\frac{\ln(q)}{\ln(pq)}\right)\right] \left[\lim_{t \to \Delta^+} A'(t) - \lim_{t \to \Delta^-} A'(t)\right],$$

which completes the proof. □

**Remark 5** *Note that if $A'(t)$ is continuous at the point $t = \frac{\ln(q)}{\ln(pq)}$, then $\lambda^C(q|p) = 0$. Thus an EV copula has non-zero $(p,q)$-quantile dependence for the values $p$ and $q$ such that $A'(t)$ is discontinuous at the point $t = \frac{\ln(q)}{\ln(pq)}$. For $p = q$ we have $\frac{\ln(q)}{\ln(pq)} = \frac{1}{2}$ and*

$$\lambda^C(p|p) = p^{2A(\frac{1}{2})-1}[A'(\tfrac{1}{2}^+) - A'(\tfrac{1}{2}^-)].$$

**Example 6** *If $A(t) = 1 - \min(\theta, \theta(1-t))$ for $\theta \in [0,1]$, then we obtain the Cuadras-Auge family of copulas [11] given by,*

$$C_\theta(u,v) = [\min(u,v)]^\theta (uv)^{1-\theta},$$

*for every $(u,v) \in [0,1]^2$. The quantity $\Delta$ in Proposition 9 is equal to $\frac{1}{2}$ for $p = q$, and $A'(\frac{1}{2}^+) - A'(\frac{1}{2}^-) = 2\theta$. Thus for $p, q \in [0,1]$, the value of the quantile dependence coefficient for this copula is given by $\lambda^{C_\theta}(q|p) = \theta p^{1-\theta} I_{[p=q]}$. Note that for this family of copulas $\lambda_L = 0$ and $\lambda_U = \theta$.*

**Example 7** *Consider an extreme value copula $C$ with the dependence function given by $A(t) = \max(t, 1-t, \theta)$, with $\theta \in [1/2,1]$. The function $A(t)$ is continuous on (0,1) and its derivative*

$$A'(t) = \begin{cases} -1, & t < 1-\theta, \\ 0, & 1-\theta < t < \theta, \\ 1, & t > \theta, \end{cases}$$

*is discontinuous at the points $t \in \{\theta, 1-\theta\}$. Thus for $p, q \in (0,1)$, from Proposition 9 we have that*

$$\lambda^C(q|p) = \theta p^{\frac{2\theta-1}{1-\theta}} \ I_{\{p^\theta = q^{1-\theta}\}} + (1-\theta) \ I_{\{p^{1-\theta} = q^\theta\}}.$$

*For this copula,* $\lambda_L(C) = I_{[\theta = 1/2]}$ *and* $\lambda_U = 2(1-\theta)$.

## 3.4 Archimedean copulas

Archimedean copula forms an important class of copulas that are easy to construct and have good analytical properties. A bivariate Archimedean copula has the form

$$C^\phi(u,v) = \phi^{[-1]}\{\phi(u) + \phi(v)\}, \tag{3.3}$$

for some continuous, strictly decreasing, and convex generator function $\phi:[0,1] \to [0,\infty]$ such that $\phi(1) = 0$ and the pseudo-inverse function $\phi^{[-1]}$ is defined by $\phi^{[-1]}(t) = \phi^{-1}(t)$, for $0 \le t \le \phi(0)$ and $\phi^{[-1]}(t) = 0$ for $\phi(0) < t \le \infty$. We call $\phi(.)$ strict if $\phi(0) = \infty$. In that case $\phi^{[-1]} = \phi^{-1}$. We note that, since $\phi$ is convex, then one-sided derivatives of $\phi'(t^-)$ and $\phi'(t^+)$ exist in (0,1] and [0,1), respectively. For Archimedean copulas, the lower and upper tail dependence coefficients are given by [11]

$$\lambda_L = \lim_{t \to 0^+} \frac{\phi^{[-1]}(2\phi(t))}{t},$$

and

$$\lambda_U = 2 - \lim_{t \to 1^-} \frac{1 - \phi^{[-1]}(2\phi(t))}{1-t} = 2 - \lim_{t \to 0^+} \frac{1 - \phi^{[-1]}(2t)}{1 - \phi^{[-1]}(t)}.$$

Let $\phi'(1)$ and $\phi(0)$ denote the derivatives at boundary of the domain of $\phi$. If $\phi'(1) = 0$, then $\lambda_U = 2 - (\phi^{-1} \circ 2\phi)(1)$ and if $\phi'(0) = -\infty$ then, $\lambda_L = (\phi^{-1} \circ 2\phi)(0)$. If $\phi'(1) < 0$ then $\lambda_U = 0$ and if $\phi'(0) > -\infty$ then $\lambda_L = 0$. The values of $\lambda_L$ and $\lambda_U$ for Archimedean copulas could be seen in [11]; see, Example 5.22.

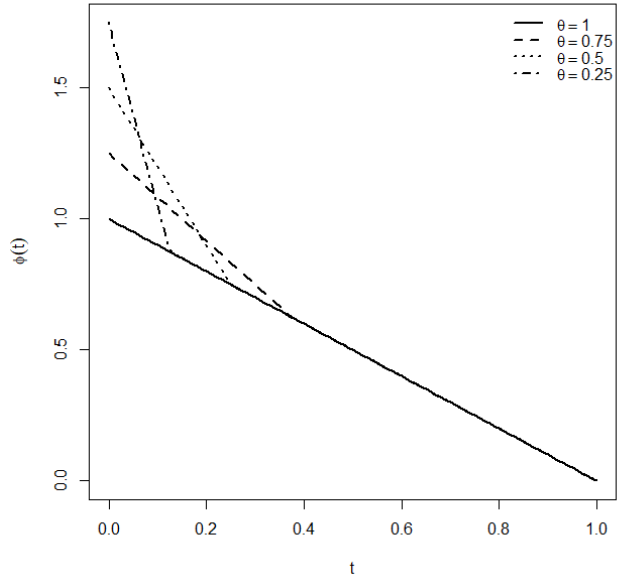

Figure 4: Plot of generator $\phi$ in Example 8.

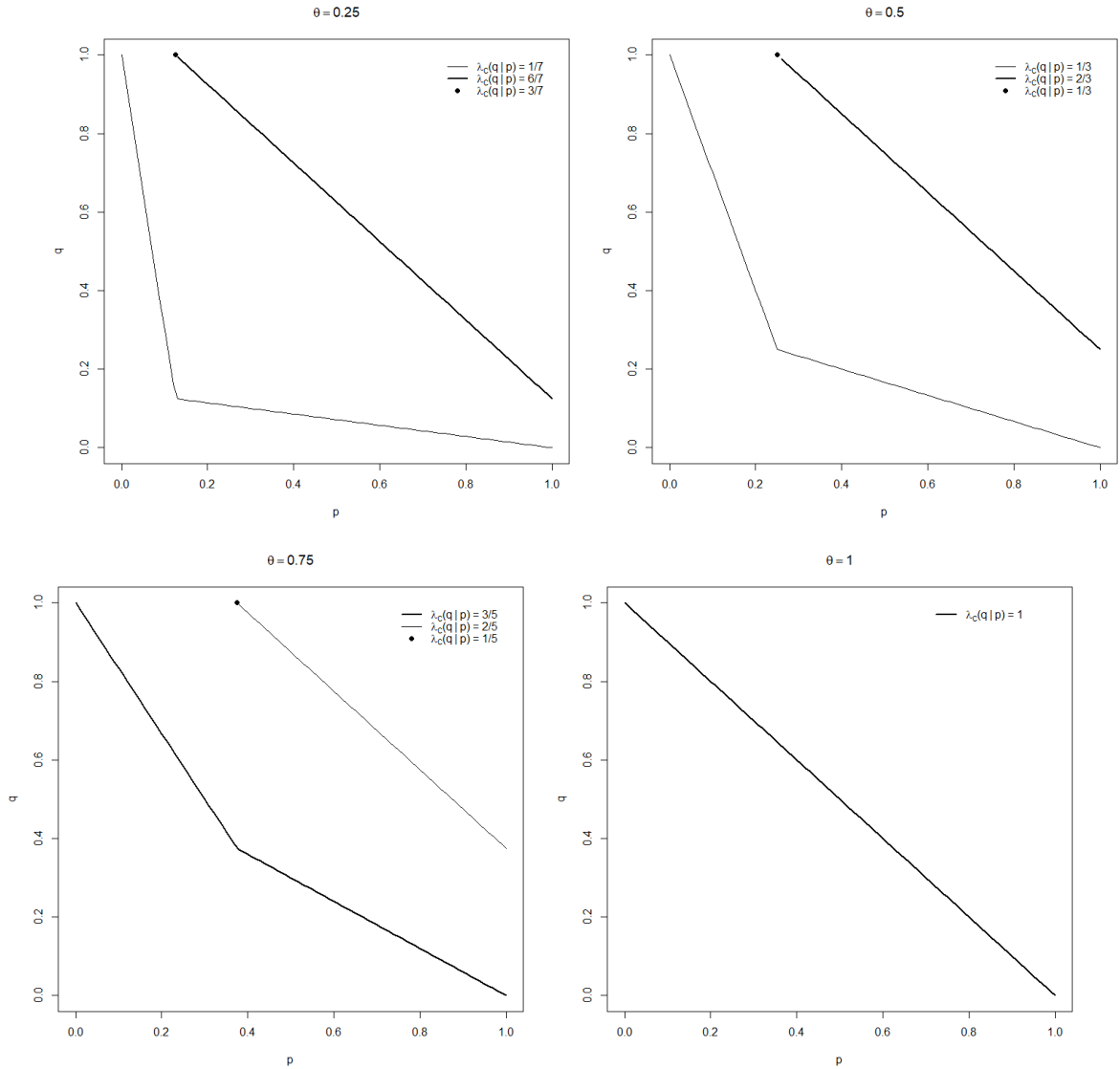

Figure 5: Plots of $\lambda^C(q|p)$ for the Archimedean copula given in Example 8.

The value of $(p, q)$-quantile dependence coefficient for Archimedean copulas can be calculated using Corollary 1 and Corollary 2. For some Archimedean copulas, the upper and lower tail dependence coefficients are equal to zero, but the quantile dependence coefficient can be non-zero. As the following example shows.

**Example 8**  *for $\theta \in (0,1]$, consider the function*

$$\phi(t) = \begin{cases} \frac{2-\theta}{\theta}(\theta - t), & 0 \leq t \leq \frac{\theta}{2}, \\ 1 - t, & \frac{\theta}{2} \leq t \leq 1, \end{cases}$$

*with the pseudo-inverse*

$$\phi^{[-1]}(t) = \begin{cases} 1-t, & 0 \le t \le 1-\frac{\theta}{2}, \\ \theta - \frac{\theta}{2-\theta}t, & 1-\frac{\theta}{2} \le t \le 2-\theta, \\ 0, & t \ge 2-\theta. \end{cases}$$

*It is easy to see that $\phi$ is the non-strict generator of an Archimedean copula given by*

$$C(u,v) = \begin{cases} u - \frac{\theta}{2-\theta}(1-v), & u \in [0,\frac{\theta}{2}], v \in [\frac{\theta}{2},1], \frac{2-\theta}{\theta}(\theta - u) + (1-v) \le 2-\theta, \\ v - \frac{\theta}{2-\theta}(1-u), & u \in [\frac{\theta}{2},1], v \in [0,\frac{\theta}{2}], \frac{2-\theta}{\theta}(\theta - v) + (1-u) \le 2-\theta, \\ \frac{\theta}{2-\theta}(u+v-\theta), & u,v \in [\frac{\theta}{2},1], u+v-1 \le \frac{\theta}{2}, \\ u+v-1, & u,v \in [\frac{\theta}{2},1], u+v-1 \ge \frac{\theta}{2}, \\ 0, & 0.w. \end{cases}$$

*This copula is lower and upper tail independent, i.e., $\lambda_L(C) = \lambda_U(C) = 0$ but the quantile dependence coefficient $\lambda^C(q|p)$ is given by*

$$\lambda^C(q|p) = \begin{cases} \frac{\theta}{2-\theta}, & q + \frac{2-\theta}{\theta}p = 1, p \in [0,\frac{\theta}{2}] \quad or \quad p + \frac{2-\theta}{\theta}q = 1, p \in [\frac{\theta}{2},1], \\ 1 - \frac{\theta}{2-\theta}, & p+q = 1+\frac{\theta}{2}, p \in (\frac{\theta}{2},1], \\ \frac{1}{2}\left(1-\frac{\theta}{2-\theta}\right), & p = \frac{\theta}{2}, q = 1, \\ 0, & o.w. \end{cases}$$

Figure 4 shows the values of $\lambda^C(q|p)$ for different values of $\theta$. Note that for $\theta = 1$, this copula reduces to $W$ and $\lambda^W(q|p) = I_{[p=1-q]}$.

## 4    Conclusion

Following measuring the degree of tail dependence between two random variables in the lower-left corner and upper-right corner of their copula domain, we developed a copula-based concept of quantile dependence between two random variables to measure the degree of dependence in specific regions of the domain. The lower and upper tail dependence coefficients are special cases of this measure. The properties of the proposed quantile dependence coefficient were studied. For illustration, in some examples of copulas that have quantile dependence, the value of the proposed measure was calculated. Expressions were obtained for the calculation of the quantile dependence coefficient in the family of Gaussian copula, Student T copula, Archimedean copulas, and Extreme value copulas. Estimating the proposed quantile dependence coefficient and its applications in dependence modeling is the subject of future research. To model datasets with different dependence patterns, copulas with different dependence structures are needed. To model the data in which there is quantile dependence, copulas that have quantile dependence are needed. One line of research could be the construction of copulas with this kind of dependence, such as the copula given by (2.3).